\documentclass[11pt,a4paper]{article}
\usepackage[margin=1in]{geometry}
\usepackage{latexsym, amsfonts, amsmath}
\newcommand{\be}{\begin{equation}}
\newcommand{\ee}{\end{equation}}
\newcommand{\bea}{\begin{eqnarray}}
\newcommand{\eea}{\end{eqnarray}}
\newcommand{\bean}{\begin{eqnarray*}} 
\newcommand{\eean}{\end{eqnarray*}}

\newcommand{\brray}{\begin{array}}
\newcommand{\erray}{\end{array}}
\newcommand{\ben}{\begin{equation}{nonumber}}
\newcommand{\een}{\end{equation}{nonumber}}

\newtheorem{dfn}{Definition}[section]
\newtheorem{thm}[dfn]{Theorem}
\newtheorem{lmma}[dfn]{Lemma}
\newtheorem{ppsn}[dfn]{Proposition}
\newtheorem{crlre}[dfn]{Corollary}
\newtheorem{xmpl}[dfn]{Example}
\newtheorem{rmrk}[dfn]{Remark}
\newcommand{\bdfn}{\begin{dfn}}
\newcommand{\bthm}{\begin{thm}}
\newcommand{\blmma}{\begin{lmma}}
\newcommand{\bppsn}{\begin{ppsn}}
\newcommand{\bcrlre}{\begin{crlre}}
\newcommand{\bxmpl}{\begin{xmpl}}
\newcommand{\brmrk}{\begin{rmrk}}
\newcommand{\edfn}{\end{dfn}}
\newcommand{\ethm}{\end{thm}}
\newcommand{\elmma}{\end{lmma}}
\newcommand{\eppsn}{\end{ppsn}}
\newcommand{\ecrlre}{\end{crlre}}
\newcommand{\exmpl}{\end{xmpl}}
\newcommand{\ermrk}{\end{rmrk}}

\newcommand{\IC}{\mathbb{C}}

\newcommand{\IE}{{I\! \! E}}

\newcommand{\IR}{\mathbb{R}}

\newcommand{\cla}{{\cal A}}
\newcommand{\clb}{{\cal B}}
\newcommand{\clc}{{\cal C}}

\newcommand{\cle}{{\cal E}}
\newcommand{\clf}{{\cal F}}
\newcommand{\clg}{{\cal G}}
\newcommand{\clh}{{\cal H}}
\newcommand{\cli}{{\cal I}}

\newcommand{\cll}{{\cal L}}

\newcommand{\clq}{{\cal Q}}

\newcommand{\cls}{{\cal S}}

\newcommand{\clv}{{\cal V}}

\def\a*{{\cal A}_{h,*}}
\def\B{{\cal B}(h)}
\def\B1{{\cal B}_1(h)}
\def\b{{\cal B}^{\rm s.a.}(h)}
\def\b1{{\cal B}^{\rm s.a.}_1(h)}

\newcommand{\ot}{\otimes}

\newcommand{\raro}{\rightarrow}

\newcommand{\lgl}{\langle}
\newcommand{\rgl}{\rangle}

\def \qed {$\Box$}

\def\a*{{\cal A}_{h,*}}
\def\B{{\cal B}(h)}
\def\B1{{\cal B}_1(h)}
\def\b{{\cal B}^{\rm s.a.}(h)}
\def\b1{{\cal B}^{\rm s.a.}_1(h)}

\begin{document}
\begin{center}
{\Large{\bf Non-existence of genuine (compact) quantum symmetries  of compact, connected smooth 
manifolds }}\\ 
{\large {\bf Debashish Goswami\footnote{Partially supported by J C Bose Fellowship from D.S.T. (Govt. of India).} }}\\ 
Indian Statistical Institute\\
203, B. T. Road, Kolkata 700108\\
Email: goswamid@isical.ac.in,\\ 
Phone: 0091 33 25753420, Fax: 0091 33 25773071\\
\end{center}
\begin{abstract}
Suppose that a compact quantum group 
${\mathcal Q}$ acts faithfully  on a smooth, compact, connected  manifold
$M$, i.e. has a $C^{\ast}$ (co)-action $\alpha$ on $C(M)$, such that $\alpha(C^\infty(M)) \subseteq C^\infty(M, {\mathcal Q})$ and 
 the linear span of 
$\alpha(C^\infty(M))(1 \otimes {\mathcal Q})$ is dense in $C^\infty(M, {\mathcal Q})$ with respect to the Fr\'echet topology. It was conjectured by the author quite a few years ago that 
    ${\mathcal Q}$ must be commutative as a $C^{\ast}$ algebra i.e.
${\mathcal Q} \cong C(G)$ for some compact group $G$ acting smoothly on $M$.  The goal of this paper is to prove the truth of this conjecture. A remarkable aspect of the proof is the use of probabilistic techniques involving Brownian stopping time. 
 
\end{abstract}

  
  \section{Introduction}
  In this article, we settle a conjecture about quantum group actions on classical spaces, which was made by the author  in \cite{Rigidity} quite a few years ago and which has been proved 
   in certain cases by him and others over the recent years. Let us give some background before stating it. 
 
 Quantum groups have their origin in both physics and mathematics, as  generalized symmetry objects of possibly noncommutative spaces.
 Following pioneering works by Drinfeld \cite{drinfeld} Jimbo \cite{jimbo}, Faddeev-Reshetikhin-Takhtajan \cite{frt} and others (see, e.g. \cite{soi}) 
  in the algebraic framework and later Woronowicz \cite{Pseudogroup}, Podles \cite{Podles}, Vaes-Kustermans \cite{kustermans} and others in the analytic setting, 
 there is by now a huge and impressive literature on quantum groups. In   \cite{manin_book}, Manin studied quantum symmetry in terms of certain universal Hopf algebras. In the 
 analytic framework of compact quantum groups a la Woronowicz, Wang, Banica, Bichon, Collins (see, e.g. \cite{ban_1}, \cite{bichon}, \cite{wang}) and many other mathematicians formulated 
  and studied quantum analogues of permutation and automorphism groups for finite sets, graphs, matrix algebras etc. This motivated the more recent theory of quantum isometry groups \cite{Goswami}
   by the 
   author of the present article in the context of Connes' noncommutative geometry (c.f. \cite{con}), which was developed further by many others  
    including Bhowmick,  Skalski, Banica, Soltan, De-Commer, Thibault, just to name a few (see, e.g.  \cite{skalski_bhow}, \cite{Soltan} etc. and the references therein).

    In this context, 
    it is important to study quantum symmetries of classical spaces. One may hope that there are many more genuine quantum symmetries of a given classical space than classical group
     symmetries which will help one understand the space better. By `genuine' we mean that the underlying algebra structure of the quantum group is noncommutative. In this context, one may mention Wang's discovery of infinite dimensional quantum permutation group $\cls^+_n$ of 
      a finite set with $n$  points where $n\geq 4$  and  the discussion on `hidden symmetry in algebraic geometry' in Chapter 13 of \cite{manin_book}.
      It follows from Wang's work that any disconnected  space with $4$ or more homoemorphic components will admit a faithful quantum symmetry given 
       by a suitable quantum permutation group. It is more interesting to look for  nontrivial and 
       interesting examples of (faithful) (co)-actions of genuine quantum groups on connected classical topological spaces as well as connected algebraic varieties. 
       Indeed, several such examples are known by now, which include:\vspace{2mm}\\
       (i) (co)action of $\cls^+_n$ on the connected compact space formed by topologically gluing $n$ copies of 
        a given compact connected space  \cite {Huchi};\\
        (ii) Co-action of the group $C^*$ algebra $C^*(S_3)$ of  the group of permutations of $3$ objects on the coordinate ring of the variety $\{ xy=0\}$ as in \cite{Etingof_walton_1}\\
        (iii)   algebraic co-action of Hopf-algebras corresponding to genuine non-compact 
  quantum groups on commutative domains associated with affine varieties as in  \cite{Wal_Wang}(Example 2.20).\\
  (iv) Algebraic co-action of (non-commutative) Hopf algebras 
   on the coordinate ring of singular curves like the cusp and the nodal cubic, given   by Krahmer and his collaborators (see \cite{kram1}, \cite{kram2}).

  However,  one striking observation is that in each of the above examples, either the underlying space 
   is not a smooth manifold ((i), (ii), (iv)) or the quantum group is not of compact type (in (iii) and (iv)). There   seems to be a natural obstacle to
   construct genuine compact quantum group action on a compact connected smooth manifold, at least when the action is assumed to be smooth in a natural sense.   
  Motivated by the fact that
  a topological action $\beta$ of a compact group $G$ on a  smooth manifold $M$ is smooth in the sense that each $\beta_g$ is a smooth map (diffeomorphism)
   if and only if it is 
    isometric w.r.t. some Riemannian structure on the manifold, the first author of this paper and some of his collaborators and students 
     tried to compute quantum isometry groups for several classical (compact)  Riemannian manifolds including the spheres and the tori. 
   Quite remarkably,  in each of  these cases, the quantum isometry group turned out to be the same as $C(G)$ where $G$ is the corresponding isometry group.  
On the other hand,  Banica et al (\cite{banica_jyotish}) ruled out the possibility of (faithful) isometric actions of a 
 large class of compact quantum groups including $\cls_n^+$ on a connected compact Riemannian manifold. All these led the first author of the present paper to make the 
following conjecture in \cite{Rigidity}, where he also gave some supporting evidence to this conjecture considering certain class of homogeneous spaces.\\   
{\bf Conjecture I:} {\it It is not possible to have smooth faithful action of a 
genuine compact quantum group on $C(M)$ when $M$ is a compact connected smooth
manifold.}\\

There have been several results, both in the algebraic and analytic set-up, which point towards  the truth of this conjecture.  For example, 
it is verified in \cite{gafa}  under the additional  condition that the action is isometric in the sense of \cite{Goswami} for some 
 Riemannian metric on the manifold. In \cite{Etingof_walton_1}, Etingof and Walton  obtained
a somewhat similar result in the purely algebraic set-up by proving that there cannot be  any  finite dimensional
Hopf algebra having inner faithful action on a commutative domain. However, their proof does not seem to extend to the infinite dimension as it crucially
depends on the semisimplicity and finite dimensionality of the Hopf algebra. We should also mention the proof   by A. L. Chirvasitu (\cite{chirvasitu}) of 
  non-existence of genuine quantum isometry in the metric space set-up (see \cite{sabbe}, \cite{ban_1}, \cite{Metric} etc.)
  for the geodesic metric of a negatively curved, compact connected Riemannian manifold.
  
  In the present article, we settle the above conjecture in the affirmative. In fact, in a pre-print written with two other collaborators, 
   the author of the present paper  posted a claim of the proof of this fact on the archive quite a few years ago but it contained a
   crucial gap. The idea was to emulate the classical averaging trick for constructing a Riemannian metric for which the given smooth CQG action is isometric. 
    However, the idea did not work mainly because we could not prove that the candidate of the Laplacian associated to the averaged metric was a second order differential operator. 
     In the present article, we circumvent the difficulties using techniques of stopping time from the theory of probability. In fact, we follow the classical line 
      of proving locality of the infinitesimal generator of the heat semigroup using stopping time of Brownian motion on manifolds. 
 
  \brmrk  
    In some sense, our results indicate that one cannot possibly have a genuine `hidden quantum symmetry' in the sense of Manin (Chapter 13 of \cite{manin_book}) 
     for smooth connected varieties coming from compact type  Hopf algebras; i.e. one must look for such quantum symmetries given by Hopf algebras of non-compact type only.
      From a physical point of view, it follows that for a classical mechanical  system with phase-space modeled by a compact connected manifold,
 the generalized notion of symmetries in terms of (compact) quantum groups coincides with the conventional notion, i.e. symmetries coming from group actions. \ermrk

  \section{Preliminaries}
 
  \subsection{Notational convention} We will mostly follow the notation and terminology of \cite{gafa}, some of which we briefly recall here. 
  All the
Hilbert spaces are over $\mathbb{C}$ unless mentioned otherwise. 
 For a complex $\ast$-algebra $\clc$, let $\clc_{s.a.}=\{c\in\clc: c^{\ast}=c\}$. We shall denote the $C^{\ast}$
algebra of bounded operators on a Hilbert space $\clh$ by $\clb(\clh)$ and the
$C^{\ast}$ algebra
of compact operators on $\clh$ by $\clb_{0}(\clh)$. $Sp$, $\overline{Sp}$ 
stand for the
linear span and closed linear span of elements of a vector space respectively, whereas ${\rm Im}(A)$ denotes the image of a linear map. 

We will deviate from the convention of \cite{gafa} in one context: we'll use the same symbol $\ot$ for any kind of topological
 tensor product, namely minimal $C^*$ tensor product, projective tensor product of locally convex spaces as well as 
  tensor product of Hilbert spaces and Hilbert modules. However, $\ot_{\rm alg}$ will be used for algebraic tensor product of vector spaces, 
   algebras or modules. 
A scalar valued inner product of Hilbert spaces will be denoted by $<\cdot, \cdot >$ and some (non-scalar)  $\ast$-algebra 
valued inner product of Hilbert modules over locally convex $\ast$-algebras  will be denoted  by $\lgl\lgl \cdot, \cdot \rgl\rgl$.
For a Hilbert $\cla$-module $E$ where $\cla$ is a $C^*$ algebra, we denote the $C^*$-algebra of adjointable right $\cla$-linear maps 
 by $\cll(E)$. In particular, we'll consider the trivial Hilbert modules of the form $\clh \ot \cla$.

   Throughout the paper, let $M$ be a compact smooth manifold of 
    dimension $m$. Let us 
   also fix an embedding of $M$ in some $\IR^n$ and let $x_1, \ldots, x_n$ denote the restriction of the canonical coordinate functions of $\IR^n$ to $M$. 
  
  \subsection{Compact quantum groups and their actions}
  We recall from \cite{gafa} and the references therein, including  \cite{Van}, \cite{Pseudogroup}, some basic facts about compact quantum groups
   and their actions. 
A compact quantum group (CQG for short) is a  unital $C^{\ast}$ algebra $\clq$ with a
coassociative coproduct 
(see \cite{Van}) $\Delta$ from $\clq$ to $\clq \ot \clq$  
such that each of the linear spans of $\Delta(\clq)(\clq\ot 1)$ and that
of $\Delta(\clq)(1\ot \clq)$ is norm-dense in $\clq \ot \clq$. 
From this condition, one can obtain a canonical dense unital $\ast$-subalgebra
$\clq_0$ of $\clq$ on which linear maps $\kappa$ and 
$\epsilon$ (called the antipode and the counit respectively) are defined, making the above subalgebra a Hopf $\ast$-algebra.

 It is known that there is a unique state $h$ on a CQG $\clq$ (called the Haar
state) which is bi invariant in the sense that $({\rm id} \ot h)\circ
\Delta(a)=(h \ot {\rm id}) \circ \Delta(a)=h(a)1$ for all $a\in\clq$. The Haar state
need not be faithful in general, though it is always faithful on $\clq_0$ at
least. The image of $\clq$ in the GNS representation of $h$ in the GNS Hilbert space $L^2(\clq,h)$ is denoted by 
 $\clq_{r}$ and it is called the reduced CQG corresponding to $\clq$. 
 
 A unitary representation of a CQG $(\clq,\Delta)$ on a Hilbert
space $\clh$ is a unitary $U \in \cll(\clh \ot \clq)$ such that the $\mathbb{C}$-linear map $V$ from $\clh$ to the Hilbert module
$\clh \ot
\clq$ given by $V(\xi)=U(\xi \ot 1)$ satisfies 
 $(V \ot {\rm id})V=({\rm id} \ot \Delta)V.$
Here, 
the map $(V \ot {\rm id})$ denotes the extension of $V \ot {\rm id}$ to the completed tensor product 
  $\clh \ot \clq$ which exists as $V$ is an isometry.

 For a Hopf algebra $H$ with the coproduct $\Delta$, 
 we write $\Delta(q)=q_{(1)}\ot q_{(2)}$ suppressing the summation notation (Sweedler's notation). For
an algebra (other than $H$ itself) or module $\cla$ and a $\mathbb{C}$-linear map $\Gamma:\cla\raro
\cla\ot_{\rm alg} H$ (typically a comodule map or a coaction) we will also use an analogue of Sweedler's notation, by writing 
$\Gamma(a)=a_{(0)}\ot a_{(1)}$. 

\bdfn
\label{CQG_action}
A unital $\ast$-homomorphism  $\alpha:\clc\raro \clc \ot \clq$, where $\clc$ is a unital $C^\ast$-algebra and 
 $\clq$ is a CQG, is
said to be an action of $\clq$ on $\clc$  if\\
1. $(\alpha \ot {\rm id})\alpha=({\rm id} \ot \Delta)\alpha$ (co-associativity).\\
2. $Sp \ \alpha(\clc)(1\ot \clq)$ is norm-dense in $\clc \ot \clq$. 
\edfn
 Given an action $\alpha$ of a CQG $\clq$ on $\clc$, there exists a norm-dense unital $\ast$-subalgebra  of $\clc$ 
 over which $\alpha$ restricts to an algebraic co-action of the Hopf algebra $\clq_0$. 
 
 An action $\alpha$  of $\clq$ on $\clc$ induces an action (say $\alpha_r$)  of the corresponding reduced CQG $\clq_r$ and the original action is faithful if and only 
  if the action of $\clq_r$ is so. 
  We say that  the action of $\alpha$ can be implemented by unitary representation if we can find a Hilbert space $\clh$ such that $\cla \subseteq \clb(\clh)$ and a unitary representation $U$  on $\clh$ such that
  $\alpha(a)=U (a \ot 1)U^{-1}$ for all $a \in \clc$. 
  It is easy to see that any unitarily implemented action is injective. In fact, as the unitary representation $U$ of $\clq$ induces a unitary representation $U_r:=({\rm id} \ot \pi_r)(U)$ of $\clq_r$ which implements $\alpha_r$, it follows that $\alpha_r$ is injective as well. If $\clc$ is separable, we can prove the converse as follows. Let $\alpha_r$ be injective.  By separability, we can find a faithful state $\phi$ on $\clc$ and then `average' it w.r.t. the faithful Haar state $h$  of $\clq_r$, i.e. define $\overline{\phi}=(\phi \ot h) \circ \alpha_r$, which is clearly $\clq_r$-invariant and also faithful. By invariance, the map $a \ot q \mapsto \alpha_r(a)(1 \ot q)$, $a \in \clc,~q \in \clq_r$ extends to a unitary representation of $\clq_r$ on the GNS space $L^2(\clc,\overline{\phi})$. This unitary representation  implements $\alpha_r$. 
  
  In the special case $\clc=C(X)$ where $X$ is a compact Hausdorff space, the above invariant state will correspond to a faithful Borel measure, say $\mu$, so that the injective reduced action is implemented by a unitary representation in $L^2(X,\mu)$.

  \subsection{ Sesquilinear form associated to a nondegenerate, conditionally positive definite, local operator}
  \bdfn
  \label{local_etc}
  Consider a linear map $\cll$ from $C^\infty(M)$ to $C(M)$ satisfying $\cll(1)=0$. We say that $\cll$ is  \\
  (i) real, if $\cll(\overline{f})=\overline{\cll(f)}$ for all $f \in C^\infty(M)$;\\
 (ii) local, if for any $x \in M$ and any $f \in C^\infty(M)$ such that $f(y)=0$ for all $y$ in an open neighbourhood of $x$, we must have $\cll(f)(x)=0$;\\
 (iii) conditionally positive definite if $\cll$ is real and  for any $f_1, \ldots, f_k \in C^\infty(M)$, $k \geq 1$ and $x \in M$, 
 the $k \times k$ matrix $(( k_\cll(f_i,f_j)(x) ))$ is 
  nonnegative definite, where $k_\cll(f,g):=\cll(\overline{f}g)-\cll(\overline{f})g-\overline{f}\cll(g).$\\
  \edfn
  
  The following result is  perhaps well-known, but we give a complete proof as we could not locate a precise reference of the result stated in this form.
  \bppsn
  \label{nondeg_loc}
  Let $\cll$ be a  local, conditionally positive definite $\cll$ with $\cll(1)=0$ as above.  Then 
   there is  a unique $C(M)$-valued, non-negative definite  sesquilinear form $<< \cdot, \cdot >>$ on $\Omega^1(M)$ (the space of smooth one-forms) such that $<<df,dg>>=k_\cll(f,g)$ for all $f,g \in C^\infty(M)$.
   
  
  \eppsn
  {\it Proof:}\\
  It is easy to see that $k_\cll(f,g)=k_\cll(g,f)$ for $f,g $ real and also $k_\cll(f,g)(x)=0$ if $f$ (or $g$) is zero on an open neighborhood of $x$. Hence $k_\cll(f,g)(x)$ depends only on the values of $f,g$ in an open neighbourhood of $x$.
  Moreover, as $\cll(1)=0$, \be \label{30} k_\cll(f,1)=k_\cll(1, f)=0\ee  for any $f \in C^\infty(M)$.
  
  Now, fix any $x \in M$. Also, fix any positive integer $k$ and smooth real functions $f_1, \ldots, f_k$ such that $f_i(x)=0$ for each $i.$
   Consider a  linear map $\Theta :M_k(C^\infty(M) )\raro M_k(\IC)$
   given by, $\Theta(G)=\left(\left( \cll(f_if_jg_{ij})(x) \right)\right) \in M_k(\IC)$, where 
   $G=(( g_{ij})) \in M_k(C^\infty(M))$. 
   By  (iii) of Definition \ref{local_etc},
   we have the following, where  $c_1, \ldots, c_k$ are complex numbers,  $H=(( h_{ij} ))\in M_k(C^\infty(M))$ and $G=H^*H$:
   \bean \lefteqn{\sum_{ij} \overline{c_i}c_j\cll(f_if_jg_{ij})(x)}\\
   &=& \sum_{i,j,p=1}^k \overline{c_i}c_j\cll(f_if_j\overline{h_{pi}}h_{pj})(x)\\
   & & \geq \sum_{i,j,p=1}^k \overline{c_i}c_j\left( f_i(x)\overline{h_{pi}}(x)\cll(f_jh_{pj})(x)+\cll(\overline{h_{pi}}f_i)(x)f_j(x)h_{pj}(x)\right)=0,\eean  
   which proves that $\Theta$ 
   is a positive linear map. 
    Note that  $M_k(C^\infty(M))$ is a $\ast$-subalgera of $M_k(C(M))$  which is unital and closed under the holomorphic functional calculus. Hence  the 
  positive linear map $\Theta$ on $M_k(C^\infty(M))$ extends uniquely  to $M_k(C(M))$ as a positive linear map, denoted again by $\Theta$. As $M_k(C(M))$ is a unital $C^*$ algebra, the (extended) positive map $\Theta$ 
  is norm-bounded with the norm $\| \Theta \|=\| \Theta(1)\|$. 
   This gives us $\| (( \cll(f_i f_jg_{ij})(x) ))\|\leq \| ((g_{ij} )) \|_\infty \|(( \cll(f_if_j)(x) ))\|$ for all $(( g_{ij}))
  \in M_k(C^\infty(M))$. But as $\cll$ is local, $\cll(f_if_jg_{ij})(x)$ depends only on the values of $ f_if_jg_{ij}$
   in an arbitrarily small open neighbourhood of $x$. If   $g_{ij}(x)=0$ for all $i,j$, then 
    for any $\epsilon>0$ we can choose  open neighbourhoods $V,W$ (say) of $x$ such that $\overline{V} \subset W$ and 
    $\| (( g_{ij}(y)))\| \leq \epsilon$ for all $y \in W$. Let $\chi$ be a smooth function supported in $W$ with $0 \leq \chi \leq 1$ and 
     $\chi|_V\equiv 1$. 
     Let $G_1(y)\equiv (( g^1_{ij}(y) ))=\chi(y) (( g_{ij}(y) )) \forall y \in M$. It  satisfies $G_1(y)=(( g_{ij}(y)))$ for all $y\in V$ and $\| G_1(y) \| \leq \epsilon$ 
     for all $y \in M$. Thus, $\| \Theta(G) \|=\| \Theta(G_1)\| \leq \epsilon \|\Theta(1)\| $, 
      i.e.$\| (( \cll(f_if_jg_{ij})(x))) \| \leq  \epsilon \| ((  \cll(f_if_j)(x) ))\|.$
    As $\epsilon$ is arbitrary, this proves $\cll(f_if_jg_{ij})(x)=0$. It also follows that \be \label{*} k_\cll(f_ig_i,f_jg_j)(x)=0\ee if 
    $f_i,g_i$ real valued smooth functions with $f_i(x)=g_i(x)=0$. 
    
    Next,  choose a local coordinate $(U, \xi_1, \ldots, \xi_m)$ (say) around $x$. 
    Without loss of generality, we can assume $\xi_i(x)=0$ for each $i$, because $k_\cll(f,g)(x)=k_\cll(f-f(x)1,g-g(x)1)(x)$ for any $f,g \in C^\infty(M)$. Choose another open neighbourhood $V_1$ of $x$ such that 
     $\overline{V_1}\subset U$ and a smooth positive function $\chi$ supported in $U$ such that $\chi|_{\overline{V_1}}\equiv 1$. 
     Now, given a real valued smooth $f$ we can write $f=f(x)1+\sum_i \partial_i(f)(x)\xi_i+R_f$ on $U$, where $\partial_i(f)(x)$ denotes the partial 
     derivative if $f$ w.r.t. the coordinate $\xi_i$ at $x$ and $R_f$ is defined in $U$. 
     Using the local Taylor expansion of $f$ around $x$ we can write $R_f=\sum_i \xi_i h_i$ where 
      $h_i$ are smooth functions defined on $U$ with $h_i(x)=0$.  
      Writing $\tilde{\phi}=\chi \phi$ for any smooth function defined at least in $U$ ( so that $\tilde{\phi} \in C^\infty(M)$), we get  
      $\tilde{f}=f(x)1+\sum_i \partial_i(f)(x)\tilde{\xi}_i+\tilde{R_f}.$ 
      As $\tilde{f}=f$, $\tilde{g}=g$ on $V_1$, we have $k_\cll(f,g)(x)=k_\cll(\tilde{f}, \tilde{g})(x)$. It also follows from (\ref{*}) 
             that $k_\cll(\tilde{h_i} \tilde{\xi}_i, \tilde{h_j} \tilde{\xi}_j)(x)=0$, hence also 
        $k_\cll(\tilde{R_f}, \tilde{R_f})(x)=0$.   By positive definiteness of $k_\cll$, 
       we have $|k_\cll(\phi, \tilde{R_f})(x)|^2\leq k_\cll(\phi, \phi)(x)k_\cll(\tilde{R_f},\tilde{R_f})(x)=0.$
             Using this as well as 
       (\ref{30}),  we get 
     \be \label{1}k_\cll(f,g)(x)=\sum_{i,j} \partial_i(f)(x)\partial_j(g)(x) k_\cll(\tilde{\xi_i},\tilde{\xi_j})(x).\ee
     
      Define a real-valued, non-negative definite  bilinear form $(\cdot,\cdot)_x$ on the cotangent space at $x$ by   setting $$ (d \tilde{\xi}_i|_x,d\tilde{\xi}_j|_x)_x=g_{ij}(x)$$ on the basis  
      $\{ d \tilde{\xi_i}|_x, i=1, \ldots, m\}$, where $g_{ij}=k_\cll(\tilde{\xi}_i, \tilde{\xi}_j)$.  
      To see the well-definedness, i.e. independence on the choice of coordinates, it suffices to 
       note that for another set of local coordinates $(\eta_1, \ldots, \eta_m)$ around $x$, we have $$d\eta_i|_x=\sum_{j=1}^m \frac{\partial \eta_i}{\partial \xi_k}(x) d\xi_k|_x,$$ and then 
         (\ref{1}) implies 
       the following: 
       $$k_\cll(\tilde{\eta}_i, \tilde{\eta}_j)(x)=
       \sum_{kl} \frac{\partial \eta_i}{\partial \xi_k}(x)
       \frac{\partial \eta_j}{\partial \xi_l}(x)g_{kl}(x).$$ 
       We can complexify $( \cdot, \cdot)_x$ to get a complex  valued sesquilinear form $<\cdot, \cdot>_x$  on the complexified cotangent spaces and using it, define a $C(M)$-valued non-negative definite  sesquilinear form given by $<< \omega, \eta>>(x)=<\omega(x), \eta(x)>_x$ for all $x \in M$. It is clear from the definition 
        that $<< df, dg>>=k_\cll(f,g)$ for all $f,g \in C^\infty(M)$. 
       \qed\\
       
         
  
  \brmrk
  \label{riem_nondeg}
  Suppose that  the map $\cll$ in the statement of Proposition \ref{nondeg_loc} also satisfies the following point-wise nondegeneracy and smoothness condition: for any $p \in M$, there are  smooth functions $f_1, \ldots, f_m$ which give a set local coordinates around $p$, $\cll(f_i)$ is smooth for each $i$ (at least around $p$)  and $(( k_\cll(f_i,f_j)(p)))$ 
   is an invertible $m\times m$ matrix. Then it is easy to see from the proof of Proposition \ref{nondeg_loc} that  there is a Riemannian structure on $M$ such that the sesquilinear form constructed in this proposition is the inner product corresponding to this Riemannian metric.
  \ermrk
   
  \subsection{Martingales and Brownian flows on manifolds}
  We will need some standard results about the Brownian motion on a compact Riemannian manifold which we briefly summarize here. For the definition, construction and 
   properties of this stochastic process, we refer to \cite{elworthy}, \cite{hsu}, \cite{protter} and the references therein. Let us consider the Riemannian 
    structure on $M$ inherited from the Euclidean Riemannian structure of $\IR^n$ and  follow the construction of \cite{hsu}, page 11, Subsect. 1.4, namely define 
     $X_t$ to be the unique solution of the stochastic differential equation $dX_t=\sum_{i=1}^n X_t P_i(X_t) \circ dW_i(t),~~X_0 \in M$, in the notation of \cite{hsu}. Here, 
      $P_i(x)$ denotes  the projection of the $i$-th coordinate unit vector of $\IR^n$  on the
tangent space $T_x M $ and $(W_1(t), \ldots, W_n(t))$ denotes the standard Brownian motion of $\IR^n$ starting at the origin. 
In this picture, $X_t$ is a process on the sample space $(\Omega, \clf, P)$ (say) of the standard $n$-dimensional Brownian motion. 
   Let $X_t(x, \omega)$ be the process `starting at $x$', i.e. the solution with $X_0=x$. 
    Let $\cll=\sum_i P_i^2$ be the Laplacian on $M$. It is known that  the Markov semigroup (`heat semigroup') given by  $T_t(f)(x):=\IE_P(f(X_t(x, \cdot)))$ has 
     $\cll$ as the infinitesimal generator. 
   We also need the following fact, which can be seen from \cite{elworthy}, Prop. 4C, Chap. I:
    \bppsn
    \label{BM_flow}
      For almost all $\omega$ in the sample space, the following hold:\\
   (i) The random  map $\gamma_t(\omega)$ given by $x \mapsto X_t(x, \omega)$ is a diffeomorphism 
     for every $t$,\\
     (ii) $(x,t) \mapsto X_t(x, \omega)$ is continuous.\\
     (iii) $X_{t+s}(x,\omega)=X_t(X_s(x,\omega),\omega).$
     \eppsn

      Let $Z$ be a separable Banach space.   We restrict our attention to separable Banach space valued random variables to avoid measure-theoretic difficulties.
      For example, the notion of Bochner or strong measurability and weak measurability coincide for separable Banach-space valued 
       random variables. We refer the reader to the lecture note by Pisier \cite{pisier} for some more details of Banach space valued measurable functions and related 
        topics. We note the following simple but useful fact. Let  $X$ be  a $Z$-valued random variable on a probability space which is almost surely norm-bounded, i.e. there is a constant $C>0$ such that $\| X(\omega) \| \leq C$ for almost all $\omega$, and let $\phi$ be  any bounded linear functional on $Z$. Then  we have $\IE(\phi(X))=\phi (\IE(X))$.
     
     We will need the concept of stopping  time (or stop time) and a  version of Doob's Optional Sampling Theorem  suitable
      for us. 
     Let us briefly 
      recall here the basics (see also \cite{protter} and the references therein). 
       We assume the usual hypotheses such as the right continuity of the filtrations considered. 
      \bdfn
      \label{stop}
      A stopping time adapted to a filtered probability space $(\Sigma, \clg, (\clg_t)_{t \geq 0}, P)$  
       is a random variable $   \tau  : \Sigma \raro \IR_+$  satisfying 
      $ \{\omega \in \Sigma ~:~\tau (\omega )\leq t\}\in \clg_t$
      for all $t \in \IR_+$.

       A family $(M_t)_{t \geq 0}$ of $Z$-valued random variables on the above filtered probability 
       space is called a $(\clg_t)$-martingale (simply martingale if the filtration is understood) 
       if $\IE(\| M_t\|)<\infty$ for each $t$, $M_t$ is  adapted to $\clg_t$ in the sense that $M_t$ is measurable w.r.t. 
        $(\Sigma, \clg_t)$ and $\IE(M_t|\clg_s)=M_s$ (almost surely) for all $0 \leq s \leq t <\infty$, where $\IE(\cdot|\clg_s)$ denotes the 
        conditional expectation with respect to $\clg_s$.

      \edfn
    
      Clearly, a Banach space valued family of random variable $M_t$ is martingale if and only if 
      for every bounded linear functional $\phi$ on $Z$, the complex valued process $\phi(M_t)$ is a martingale in the usual classical sense. 
      Adapting the proof of the classical Optional Sampling Theorem, we get the following version of Theorem 18 of Chapter I, page 10 of \cite{protter}:
     
     \bppsn
   
     \label{optsampl}
     Let $(M_t)$ be a $Z$-valued right continuous (i.e. for almost all $\omega$, $t \mapsto M_t(\omega)$ is right continuous) martingale as above. 
     Then for any bounded stopping time, the process $M_{\tau \wedge t}$ is a martingale, 
       where $a \wedge b:={\rm min}(a,b).$
     \eppsn
  {\it Proof:}\\ Let $t_0>0$ be some constant such that  $\tau \leq t_0$  almost surely. It is enough to prove that $\phi (M_{\tau \wedge t})$ is a scalar-valued martingale for each bounded linear functional $\phi$ on $Z$. But this follows by
   applying Theorem 18 of Chapter I of \cite{protter} to the scalar-valued martingale $\phi(M_{t \wedge t_0})$, or applying 
    Problem 3.23 (part (i)) of Chapter 1, page 20 of \cite{karatzas} to $\phi(M_t)$. \qed\\

  \section{Main results}
   Throughout this section, let  $M$ be a compact smooth manifold of dimension $m$ and $\clq$ be a CQG with a faithful action $\alpha$ on $C(M)$.
   
  \subsection{Smooth action}
  
 We refer to \cite{gafa} for a detailed discussion on the natural Fr\'echet topology 
  of $C^\infty(M)$ as well as the space of $\clb$-valued smooth functions $C^\infty(M, \clb)$ for any Banach space $\clb$. Indeed, by the nuclearity 
   of $C^\infty(M)$ as a locally convex space, $C^\infty(M, \clb)$ is the unique topological tensor product of $C^\infty(M)$ and $\clb$ in the category
    of locally convex spaces. This allows us to define $T \ot {\rm id}$ from $C^\infty(M, \clb)$ for any Fr\'echet continuous linear 
     map $T $ from $ C^\infty(M) $ to $C^\infty(M)$ (or, more generally, to some other locally convex space).   We also recall from \cite{gafa} the space $\Omega^1(M)\equiv 
      \Omega^1(C^\infty(M))$ of smooth one-forms and the space $\Omega^1(M, \clb)$ of smooth $\clb$-valued one-forms, as well as 
       the natural extension of the differential map $d$ to a Fr\'echet continuous map from $C^\infty(M, \clb)$ to $\Omega^1(M, \clb)$. In fact, for 
        $F \in C^\infty(M, \clb)$, the element $dF \in \Omega^1(M, \clb)$ is the unique element satisfying $({\rm id } \ot \xi)(dF(m))=(dF_\xi)(m),$ 
         for every continuous linear functional $\xi$ on $\clb$, where $m \in M,$ $dF(m) \in T^*_m M \ot_{\rm alg} \clb$ and $F_\xi \in C^\infty(M)$ is given by $F_\xi(x):=\xi(F(x))~\forall x \in M.$

 We now define a smooth action following \cite{gafa}.
\bdfn 
 In case $\clc=C(M)$, where $M$ is a smooth compact manifold, we say that an action $\alpha$ of a CQG $\clq$ on $C(M)$ is smooth if 
 $\alpha$ maps $C^\infty(M)$ into $C^\infty(M, \clq)$ and   $Sp \ \alpha(\clc)(1\ot \clq)$ is dense in $C^\infty(M, \clq)$ in the Fr\'echet topology.
 We say that the action is faithful if the algebra generated by elements of the form 
 $\alpha(f)(p)\equiv ({\rm ev}_p \ot {\rm id})(\alpha(f))$, where $f \in C(M), p \in M$ is norm-dense in $\clq$.
 \edfn
 \brmrk
 In case $\clq=C(G)$ where $G$ is a compact group  acting on $M$, say by  $\alpha_g: x \mapsto gx$, the  smoothness of the induced action $\alpha$ given by $\alpha(f)(x,g)=f(gx)$ on $C(M)$ in the sense of the above definition means the smoothness of the map $M \ni x \mapsto gx$  for each $g$.
 \ermrk

 It has been proved in \cite{gafa}, following arguments of \cite{Podles}, \cite{free} etc. 
  that given any smooth action $\alpha$ of $\clq$ on $C(M)$ there is a Fr\'echet dense unital $\ast$-subalgebra
 $\clc_0$ of $C^\infty(M)$ on which $\alpha$ restricts to an algebraic co-action of $\clq_0$. 
 It also follows (see \cite{injective}, Corollary 3.3) that  for any smooth action $\alpha$, the corresponding reduced action $\alpha_r$ is injective and hence it is  implemented by some unitary representation.

Suppose that $M$ has a Riemannian structure with the corresponding  $C^{\infty}(M)$-valued inner product $\lgl\lgl\cdot, \cdot\rgl\rgl$ on $\Omega^1(C^\infty(M))$  as in \cite{gafa}.
If $\cll$ is the Laplacian corresponding to 
  the Riemannian structure and $k_\cll(f,g):=
  \cll(\overline{f}g)-\cll(\overline{f})g-\overline{f}\cll(g)$ for $f,g \in C^\infty(M)$, 
   we have  $<<df, dg  >>=k_\cll(f,g)$.

  However, we need to consider more general non-negative definite sesquilinear forms  on
    $\Omega^1(M)$, possibly $C(M)$-valued ones, as already encountered and discussed in Subsection 2.3.
 For any such sesquilinear form  $<<\cdot, \cdot>>^\prime$, there is a canonical $C(M,\clq)$-valued, non-negative definite  sesquilinear form on 
  $\Omega^1(C^\infty(M)) \ot \clq$ determined by 
  $$<< \omega_1 \ot q_1, \omega_2 \ot q_2>>^\prime=
  <<\omega_1,\omega_2>>^\prime q_1^*q_2.$$ If $<<\cdot, \cdot>>^\prime$ is an inner product, i.e. nondegenerate on $\Omega^1(C^\infty(M))$ or on some smaller subspace of it, the corresponding sesquilinear form on 
  $\Omega^1(C^\infty(M),\clq)$ (or, its restriction  on a suitable subspace ) is an inner product too.

   We extend the definition of \cite{gafa} for Riemannian inner product preserving actions to a more general setting of  $C(M)$-valued,
   non-negative definite sesquilinear form , which is Fr\'echet-continuous, i.e. $(\omega, \eta )\mapsto <<\omega, \eta>>^\prime$ is a continuous map from the Fr\'echet space $\Omega^1(C^\infty(M)) \times \Omega^1(C^\infty(M))$  to $C(M)$.
   \bdfn 
A smooth action $\alpha$ on $M$ is said to preserve  a $C(M)$-valued, Fr\'echet continuous  sesquilinear form $<< \cdot, \cdot>>^\prime$  if 
\begin{eqnarray}
\label{inner_prod_pres_111}
 \lgl \lgl d\alpha (f), d\alpha (g)\rgl\rgl^\prime=\alpha(\lgl\lgl df,dg\rgl\rgl^\prime )
 \end{eqnarray}
 for all $f,g\in C^\infty(M)$. 
\edfn
 It is easy to see, in case by the  Fr\'echet continuity of the maps involved that it is enough to have  (\ref{inner_prod_pres_111}) 
 for  $f,g \in \clc_0.$ 
 
 \subsection{Averaging of the Riemannian metric}
  
Let $M$ be as  before and let $\alpha$ be a faithful smooth action of $\clq$  on $C(M)$.         
    Replacing $\clq$ by $\clq_r$ we can assume without loss of generality that $\clq$ has faithful Haar state and $\alpha=\alpha_r$. 
    It is also known (see \cite{huichi_2})  that $\clq_r$ is of Kac type, hence $h$ is tracial and $\kappa$ is norm-bounded 
     on $\clq=\clq_r$. 
  Let $\clq_0$ be the canonical dense Hopf $\ast$-algebra for $\clq$ and $\clc_0$ be a 
  Fr\'echet-dense unital $\ast$-subalgebra of $C^\infty(M)$ on which $\alpha$ is algebraic. 
 Moreover, as explained in Subection 2.2, 
  choose some faithful $\alpha$-invariant Borel measure $\mu$ on $M$ and the corresponding unitary representation 
  $U$  on $L^2(M,\mu)$ implementing $\alpha$, i.e. 
  $\alpha(f)=U (f \ot 1)U^{-1}$, where $f \in C(M)$ is viewed as a multiplication operator on $L^2(M, \mu)$. 
Let $L^2(\clq)$ be the GNS space of the Haar state $h$ and identify $\cll(\clh \ot \clq)$ (for any Hilbert space $\clh$) as a subalgebra of 
$\clb(\clh \ot L^2(\clq))$. The vector state $<1, \cdot 1>$ on $\clb(L^2(\clq))$ extends $h$ and we continue to denote it by $h$.

 Denote by $M_F $ and $M_f$ the operators of left multiplication by $F$ (respectively $f$) on the Hilbert $\clq$-module $L^2(M, \mu) \ot \clq$
  (respectively $L^2(M, \mu)$). Most often we may write simply $F$ or $f$ for $M_F$ or $M_f$ respectively 
   by making slight abuse of notation. 

 \blmma
 \label{888}
 For $F \in \clc_0 \ot_{\rm alg}\clq_0 \subset C^\infty(M, \clq)$,
  we have $$ ({\rm id} \ot h)(U^{-1}M_F U)=M_{F^\sharp},$$
  where $F^\sharp=({\rm id} \ot h)(U^{-1}(F)) \in \clc_0.$ 

 \elmma
 {\it Proof:}\\ 
 It is sufficient to prove the lemma for $F=f \ot q$, where $f \in \clc_0, q \in \clq_0$. Using Sweedler's notation and the trace property of $h$, 
   we have for $g \in \clc_0$: \bean \lefteqn{({\rm id} \ot h)(U^{-1}M_F U)g}\\
  &=& ({\rm id} \ot h)(U^{-1}M_F U(g \ot 1))=({\rm id} \ot h)(U^{-1}(fg_{(0)} \ot qg_{(1)}))\\
 & =& f_{(0)}g_{(0)(0)}h(\kappa(g_{(0)(1)})\kappa(f_{(1)})qg_{(1)})=f_{(0)}g_{(0)}h(\kappa(f_{(1)})qg_{(1)(2)}\kappa(g_{(1)(1)})))\\
 &=& f_{(0)}g_{(0)}h(\kappa(f_{(1)})q)\epsilon(g_{(1)})=f_{(0)}h(\kappa(f_{(1)})q)g=F^\sharp g. \eean\qed\\

 \bcrlre
 \label{cor_state}
 The map $F\mapsto \Psi(F):=({\rm id} \ot h)(U^{-1}M_F U)$ extends to a unital completely positive map from $C(M, \clq)$ to $C(M)$. In particular, 
 ${\rm ev}_p \circ ({\rm id} \ot h)(U^{-1} \cdot U)$
  extends to a well-defined state on $C(M)$. Moreover, 
    $\Psi$ is $\clq$-invariant in the sense that
  \be \label{psi_inv} (\Psi \ot {\rm id})\circ ({\rm id} \ot \Delta)
  =\alpha(\Psi(\cdot)).\ee 
 \ecrlre
 {\it Proof:}\\ The map is clearly norm-bounded and completely positive by the formula that defines it. 
  It also follows from Lemma \ref{888}  that it maps the dense subspace $\clc_0 \ot_{\rm alg} \clq_0$ into $C(M)$. 
   By norm-continuity, the image of the map must be contained in $C(M)$.
   
   To prove the invariance, it is enough to prove (\ref{psi_inv}) for $F=f \ot q$, where $f \in \clc_0, q \in \clq_0$. To this end, 
    note that as $\kappa^2={\rm id}$, we have $q_{(1)}\kappa(q_{(2)})=q_{(2)}\kappa(q_{(1)})=\epsilon(q)1$ for all $q \in \clq_0$. 
    Moreover, we have \be \label{hopf} h(q)1=(h \ot \kappa)(\Delta(q))=h(q_{(1)})\kappa(q_{(2)}),\ee  for all $q \in \clq_0$. Now, the left hand 
     side of (\ref{psi_inv}) for $F=f \ot q$ equals
     $ f_{(0)} \ot h(\kappa(f_{(1)})q_{(1)})q_{(2)}.$
          By (\ref{hopf}), $\kappa^2={\rm id}$  as well as the identity $\Delta \circ \kappa=\sigma \circ (\kappa \ot \kappa) \circ \Delta$ 
      where $\sigma$ denotes flip, 
      we have: \bean \lefteqn{ h(\kappa(f_{(1)})q_{(1)})q_{(2)}}\\
      &=& q_{(2)}(h \ot \kappa)(\Delta(\kappa(f_{(1)})q_{(1)}))\\
      &=& q_{(2)} h(\kappa(f_{(1)(2)})q_{(1)(1)})\kappa(\kappa(f_{(1)(1)})q_{(1)(2)})\\
     &=& q_{(2)}h(\kappa(f_{(1)(2)})q_{(1)(1)})\kappa(q_{(1)(2)})f_{(1)(1)}\\
     &=& h(\kappa(f_{(1)(2)})q_{(1)(1)}) q_{(2)}\kappa(q_{(1)(2)})f_{(1)(1)}\\
     &=& h(\kappa(f_{(1)(2)})q_{(1)})q_{(2)(2)}\kappa(q_{(2)(1)})f_{(1)(1)}~~({\rm as}~q_{(1)(1)} \ot q_{(2)} \ot q_{(1)(2)}
     =q_{(1)}\ot q_{(2)(2)} \ot q_{(2)(1)})\\
     &=& h(\kappa(f_{(1)(2)})q_{(1)}\epsilon(q_{(2)}))f_{(1)(1)}\\
     &=& f_{(1)(1)}h(\kappa(f_{(1)(2)})q).\eean

     It follows that \bean \lefteqn{f_{(0)} \ot h(\kappa(f_{(1)})q_{(1)})q_{(2)}}\\
     &=& f_{(0)} \ot f_{(1)(1)}h(\kappa(f_{(1)(2)})q)\\
     &=& f_{(0)(0)}\ot f_{(0)(1)}h(\kappa(f_{(1)})q)\\
     &=& \alpha(f_{(0)}h(\kappa(f_{(1)})q))\\
     &=& \alpha(\Psi(f \ot q)),\eean
     which is the right hand side of (\ref{psi_inv}).
     
   \qed\\  
 
 Choose and fix any Riemannian structure on $M$, for example, the one inherited from the embedding $M \subset \IR^n$ as in Subsection 2.4,  and write $\cll$ for the Laplacian on $M$ and 
let $$ \hat{\cll}(f)=({\rm id} \ot h)\left(U^{-1}((\cll \ot {\rm id})(\alpha(f))U\right)$$ for $f \in C^\infty(M)$.
 Here we have identified scalar or $\clq$-valued  functions with the corresponding left multiplication operators in
  appropriate Hilbert spaces or Hilbert modules, as understood from the context. By Lemma \ref{888} and Corollary \ref{cor_state},
  $\hat{\cll}(f) \in C(M)$. As $\alpha$ is Fr\'echet continuous, it is clear that $\hat{\cll}$ is continuous w.r.t. the  Fr\'echet topology 
   on $C^\infty(M)$ and the norm topology on $C(M)$. 
   We also observe that for $f \in \clc_0$, $\hat{\cll}(f)=((\cll \ot {\rm id})(\alpha(f))^\sharp=\Psi(\cll(f_{(0)})\ot f_{(1)}),$ so in particular, $\hat{\cll}(\clc_0) \subseteq \clc_0.$
   
 We now claim the following:
 \bthm
 \label{pre_main}
 $\hat{\cll}$ satisfies the hypotheses of Proposition \ref{nondeg_loc}. Moreover, the restriction of the sesquilinear form $k_{\hat{\cll}}$ induced by $\hat{\cll}$ on the Fr\'echet dense subspace $\Omega^1(\clc_0)$, which is the $\clc_0$-module generated by $\{ df,~f\in \clc_0\} $, is an inner product.
  \ethm
 {\it Proof:}\\
We divide the proof into several steps. As before,  we'll throughout make the identifications with functions (scalar or $C^*$ algebra valued) and operators of left multiplication by them 
  on appropriate Hilbert spaces or modules. Clearly, $\hat{\cll}(1)=0$. 
 Moreover, as $\Psi(F^*)=\Psi(F)^*$, we have $\hat{\cll}(\overline{f})=\overline{\hat{\cll}(f)}$ for all $f \in C^\infty(M)$. \vspace{2mm}\\
 {\it Step 1:} $\hat{\cll}$ is conditionally completely positive.\vspace{2mm}\\  
 It follows from the following claim:
 \be \label{345} k_{\hat{\cll}}(f,g)= \Psi(<<d\alpha(f), d\alpha(g)>>)\equiv 
 ({\rm id} \ot h)\left(U^{-1} (<< d\alpha(f), d \alpha(g)>>)U\right),\ee  
 where $<< \cdot, \cdot>>$ is the $C^\infty(M, \clq)$ valued inner product on $\Omega^1(M, \clq)$ coming from the Riemannian structure associated to $\cll$. Indeed, as $U$ is a unitary and $h$ is a positive functional, 
 (\ref{345}) will imply that $k_{\hat{\cll}}$ is nonnegative definite. 
  Note that $<<d\alpha(f), d\alpha(g)>>=(\cll \ot {\rm id})(\alpha(\overline{f}g))-(\cll \ot {\rm id})(\alpha(\overline{f}))\alpha(g)-\alpha(\overline{f})(\cll \ot {\rm id})(\alpha(g)).$
  
  To prove (\ref{345}), we first observe the following,  where   $\phi, \psi \in \clc_0 \ot_{\rm alg} \clq_0$:
    \bean \lefteqn{\hat{\cll}(\phi)\psi}\\
    &=&({\rm id} \ot h)\left(U^{-1}(\cll \ot {\rm id})(\alpha(\phi))U(\psi \ot 1) \right)\\
        &= &({\rm id} \ot h)\left(U^{-1}(\cll \ot {\rm id})(\alpha(\phi))\alpha(\psi)U\right).\eean
     By continuity of $\hat{\cll}$ and $\cll \ot {\rm id}$, the above equation extends to all $\phi,\psi \in C^\infty(M)$. 
     Taking $\phi=\overline{f}g$, $\psi=1$, we get 
     \be \label{aaa} \hat{\cll}(\overline{f}g)=({\rm id} \ot h)\left(U^{-1}(\cll \ot {\rm id})(\alpha(\overline{f}g))U\right).\ee Similarly, taking $(\phi, \psi)=(\overline{f}, g)$ as well as $(\phi, \psi)=(g, \overline{f})$, we get 
     \be \label{bbb} \hat{\cll}(\overline{f})g=
     ({\rm id} \ot h)\left(U^{-1}(\cll \ot {\rm id})(\alpha(\overline{f}))\alpha(g)U\right),\ee
     \be \label{ccc} \overline{f}\hat{\cll}(g)=({\rm id} \ot h)\left(U^{-1}\alpha(\overline{f})(\cll \ot {\rm id})(\alpha(g))U\right).\ee
      Combining (\ref{aaa}, \ref{bbb}, \ref{ccc}), we prove (\ref{345}) and hence {\it Step 1}.\vspace{2mm}\\

   {\it Step 2: Locality}\vspace{2mm}\\
  
  Let us consider the Brownian motion $(X_t)$  corresponding to the Riemannian structure given by $\cll$ and let $\gamma_t$ be the random 
   flow of automorphism as in the Proposition \ref{BM_flow}. For a Banach space $E$ let $L^\infty(\Omega, E)$ be the Banach space of 
    $E$-valued essentially bounded measurable functions, to be viewed as $E$-valued random variables. 
  Let $j_t : C(M, \clq) \raro  L^\infty(\Omega, C(M, \clq))$ be the $\ast$-homomorphism given by 
  $$j_t(F)(\omega)(x)\equiv j_t(F)(x, \omega)=F(X_t(x, \omega)).$$ That is, $j_t(F)(\omega)=F \circ \gamma_t(\omega).$ 
  Let $\cle$ be the $C^*$ subalgebra of $\cll(L^2(M, \mu) \ot \clq)$ given by $\cle=\{ U^{-1}M_F U:~F \in C(M, \clq)\}\equiv U^{-1}C(M, \clq)U$. 
  We note that $\cle$ is separable. Indeed, as $M$ is compact and $\clq$ acts faithfully on the separable $C^*$ algebra $C(M)$, $\clq$ is separable too. This implies the separability of 
   $C(M, \clq)\equiv C(M) \ot \clq$ and hence of  $\cle \cong C(M, \clq)$.  Moreover, $\cle$ contains $C(M) \ot 1=U^{-1} \alpha(C(M))U$. 
   
   Using the identification  of $C(M, \clq)$ with the  left multiplication operators, define $J_t : C(M, \clq) \raro L^\infty(\Omega, \cle)$ by $$ J_t(F)(\omega)=U^{-1} j_t(F)(\omega)U.$$ 
    Clearly, $J_t$ is a unital $\ast$-homomorphism.  We also have a natural 
     embedding $L^\infty(\Omega, \cle) \subseteq \clb(L^2(\Omega) \ot L^2(M, \mu) \ot L^2(\clq,h))$
      and in this picture, we can write $ J_t(F)=\tilde{U}^{-1}j_t(F)\tilde{U},$
      where $\tilde{U}= I_{L^2(\Omega)} \ot U$.

   Let $T_t$ be the Markov semigroup (heat semigroup) generated by $\cll$, which is given by the formula $T_t(f)(x)=\IE(f(X_t(x,\cdot)))$ for all $t \geq 0$.
   As $T_t$ is a $C_0$-semigroup of completely positive maps from $ C(M)$ to $C(M)$, we have the ampliation  $\tilde{T_t}:=T_t \ot {\rm id} : C(M, \clq) \raro C(M, \clq)$.
    In fact, we have $\tilde{T_t}(F)(x)=\IE(F(X_t(x, \cdot)))$ for $F \in C(M, \clq)$. 
     Let $\IE_s$ denote the conditional expectation w.r.t. the sub $\sigma$ algebra generated by $\{ X_u^{-1}(B), B \in \clb_\cle,~u \leq s\}$, where $\clb_\cle$ denotes 
      the Borel $\sigma$-algebra of the Banach space $\cle$. 
     We have \be \label{cocy} \IE_s \circ j_{s+t}=j_s \circ \tilde{T_t} \ee for $s,t\geq 0$, which follows from the Markov property of 
      the Brownian motion, e.g, as given by (iii) of Proposition \ref{BM_flow}. Let $\tilde{\cll}=\cll \ot {\rm id} $ on $C^\infty(M, \clq)$.
       Clearly, for all $F \in C^\infty(M, \clq)$, the following holds: 
       \be \label{010} \frac{d}{dt}\tilde{T}_t (F)=\tilde{\cll}\circ \tilde{T}_t(F)=\tilde{T}_t \circ \tilde{\cll}(F).\ee To verify this, we should at first 
        note that $\tilde{T}_t$ and $\tilde{\cll}$ commute as $T_t$ and $\cll$ do so. Furthermore, 
         we have \be \label{437} \tilde{T}_t(F)-F=\int_0^t \tilde{T}_s \circ \tilde{\cll}(F)ds.\ee We can  first verify (\ref{437}) on the Frechet dense subspace $\clc_0 \ot_{\rm alg}
         \clq_0$ and then extend it  to the whole of $C^\infty(M, \clq)$ by continuity of the maps involved. From this, (\ref{010}) follows 
          immediately.

     Next, we define a unital $\ast$-homomorphism $\Pi_t: C(M) \raro L^\infty(\Omega, \cle)$ by $$ \Pi_t(f)=J_t(\alpha(f)).$$ 
      For $f \in C^\infty(M)$, 
      define $$M_t^f=\Pi_t(f)-\int_0^t J_s(\tilde{\cll}(\alpha(f)))ds.$$  By  continuity of the Brownian flow, 
       the integrand on the right hand side is continuous in $s$ for almost all $\omega$ and hence convergent absolutely in the norm of $\cle$.
        We make the following:
       $${\rm Claim:}~ (M_t^f)_{t \geq 0}~{\rm is~a~continuous~\cle-valued~ martingale 
       ~w.r.t.~ the~ filtration ~of~ the~ Brownian~ motion}.$$ To prove this claim, first note that the continuity of $(M_t^f)$ follows from 
       the continuity of the Brownian flow w.r.t. the time parameter $t$. Thus, it is enough to prove the martingale property of $(M^f_t)$.
       To this end, note that   
        $\IE_u j_s(\tilde{\cll}(\alpha(f)))=j_{u} (\tilde{T}_{s-u} \circ \tilde{\cll}(\alpha(f)))=\frac{d}{ds}j_u(\tilde{T}_{s-u}(\alpha(f))),$ for $0 \leq u \leq s \leq t$, which follows from 
          (\ref{cocy}) and (\ref{010}). On the other hand, $\IE_u(j_s(F))=j_s(F)$ for $s \leq u$ and $F \in C^\infty(M)$ by definition of 
          the filtration. 
          Hence we have (for $u \leq t$ and  almost all $\omega \in \Omega$):
          \be \label{555} \IE_u(M_t^f)(\omega)
         = \IE_u(\Pi_t(f))(\omega)-U^{-1} \left( \int_0^u 
      j_s(\tilde{\cll}(\alpha(f)))(\omega)  ds+\int_u^t \frac{d}{ds} 
      j_u(\tilde{T}_{s-u}(\alpha(f)))(\omega) ds \right)U.\ee
      But observe that 
      \be \label{556} \int_0^u j_s(\tilde{\cll}(\alpha(f))) ds +\int_u^t \frac{d}{ds} 
      j_u(\tilde{T}_{s-u}(\alpha(f)))(\omega) ds=
      \int_0^u j_s(\tilde{\cll}(\alpha(f))) ds+ j_u(\tilde{T}_{t-u}(\alpha(f))
      -j_u(\alpha(f)).\ee As $j_s(\cdot)$ is measurable w.r.t. the $\sigma$-algebra $\sigma\left( X_v^{-1}(B),~B \in \clb_\cle,~ v \leq s\right)$, we have $\IE_s \circ j_s=j_s$ for all $s$. 
      Moreover,
       $\IE_u \circ j_t= j_u \circ \tilde{T}_{t-u}$, hence $\IE_u (\Pi_t(f))(\omega)=U^{-1}j_u \tilde{T}_{t-u}(\alpha(f))(\omega) U$.  We also have 
      $j_t\circ \tilde{T}_{t-u}=\IE_t \circ j_t=j_t$. 
      Combining the above observations with (\ref{556}) and 
      interchanging  $\IE_u$ with the integral by appropriate 
       continuity  of the maps involved, the right hand side of (\ref{555}) reduces to  
      $$  \Pi_u(f)(\omega)-U^{-1} \left( \int_0^u
      j_s(\tilde{\cll}(\alpha(f)))(\omega) ds\right)U=M_u(f)(\omega),$$ which proves the claim.


           Now, let $Y_i(t)=\Pi_t(x_i)$. Observe  that  $Y_i(0)=x_i\ot 1$. 
                             To show the  locality of $\hat{\cll}$ at a point $p=(p_1, \ldots, p_n) $ of $M \subset \IR^n$, consider 
           $f=\phi(x_1, \ldots, x_n),$ where $\phi$ is a 
            smooth real-valued function on $\IR^n$, and assume  that $f$ is zero on a neighbourhood of $p$. Choose 
             small enough $\epsilon_0>0$ such that $\phi(y_1, \ldots, y_n)=0$ 
             whenever $|y_i-p_i| \leq \epsilon_0$ for all $i$ and $y=(y_1, \ldots, y_n)\in M$. 
                         It is clear from the continuity properties of the Brownian flow (see Proposition \ref{BM_flow}) 
                          that $t \mapsto J_t(F)(\omega)$  is 
                          norm continuous for almost all $\omega$ and fixed $F \in C^\infty(M,\clq)$. Let  $F=\tilde{\cll}(\alpha(f))$ and 
                           let $\tau^{\prime \prime}_\epsilon(\omega)$ ($\epsilon>0$) be 
                           the infimum of $t \geq 0$ (which is defined to be $+\infty$ if no such $t$ exists) for which $\|J_t(F)( \omega)-J_0(F)( \omega) \| >\epsilon$. 
                           It is clearly a stopping time. 
                           Observe that, as $\Pi_t$ is a homomorphism, $\Pi_t(f)=\phi(\Pi_t(x_1), \ldots, \Pi_t(x_n))=\phi(Y_1(t), \ldots, Y_n(t))$.
               Furthermore, consider another  stopping time $\tau^\prime_\epsilon=\tau^\prime_\epsilon(\omega)$ to be the infimum of all $t \geq 0$ 
              for which $\|Y_i(t, \omega)-x_i \ot 1\| > \epsilon$ for some $ i$. Finally, let 
            $\tau_\epsilon={\rm min}(\tau^\prime_\epsilon,\tau^{\prime \prime}_\epsilon, 1),$ which is a bounded stopping time. 
            Applying Proposition \ref{optsampl} to the (continuous) martingale $M_t^f$, we 
            conclude that $M^f_{t \wedge \tau_\epsilon}$ is a martingale too, hence in particular, 
            $\IE(M^f_{\tau_\epsilon})=M_0^f=f \ot 1.$ In other words, $$ \IE(\Pi_{\tau_\epsilon}(f))-f \ot 1=\IE\left( \int_0^{\tau_\epsilon}
            \tilde{U}^{-1}j_s(\tilde{\cll}(\alpha(f)))\tilde{U} ds \right)=
           \IE\left( \int_0^{\tau_\epsilon}
           J_s(F) ds \right) .$$
            By definition of $\tau_\epsilon$ and continuity of the Brownian flow, it is clear that 
            $\| J_s(F)(\omega)-J_0(F)( \omega) \| \leq \epsilon$
             for all $s \leq \tau_\epsilon.$  Hence we have 
              $\int_0^{\tau_\epsilon} \| J_s(F)(\omega)-
             J_0(F)(\omega)\| ds \leq \tau_\epsilon(\omega) \epsilon.$  It follows that 
           \bean \lefteqn{\|  \IE( \int_0^{\tau_\epsilon} J_s(F)ds) -\IE(\tau_\epsilon)J_0(F)\|}\\
           &=& \| \IE \left( \int_0^{\tau_\epsilon} (J_s(F)-J_0(F))ds \right) \| \\
           &\leq  &   \IE \left(  \int_0^{\tau_\epsilon} \| (J_s(F)-J_0(F))\| ds \right)\\
           &  \leq &  \epsilon \IE(\tau_\epsilon),\eean hence 
            \be \label{loc_formula} 
            \lim_{\epsilon \raro 0+}\frac{\IE(\Pi_{\tau_\epsilon}(f))-f\ot 1}{\IE(\tau_\epsilon)}
            =\lim_{\epsilon \raro 0+}\frac{\IE(\int_0^{\tau_\epsilon} J_s(F)ds)}{\IE(\tau_\epsilon)}
            =\IE(J_0(F))=U^{-1}\tilde{\cll}(\alpha(f))U,\ee
             where the convergence is in the norm of $\cll(L^2(M,\mu) \ot \clq).$

             For a fixed $t$ and $\omega$, let us denote by $\clb_{t, \omega} \subseteq \cll(L^2(M, \mu) \ot \clq)$ the commutative unital $C^*$ algebra generated by $\{ \Pi_t(f)(\omega), g \ot 1,$ $f,g \in C(M)\}$.
               Clearly, $\clb_{0, \omega}=C(M) \ot 1$ is a common $C^*$-subalgebra of all $\clb_{t, \omega}$. Let $\cls$ be the (convex, weak-$\ast$ compact)
               set of states $\zeta$ on $\clb_{t, \omega}$ which extends ${\rm ev}_p$ on $C(M) \ot 1 \cong C(M)$, i.e. $\zeta(g \ot 1)={\rm ev}_p(g \ot 1)=g(p)~\forall g$. 
               By standard arguments we can prove that any extreme point of $\cls$ is 
                also an extreme point of the set of all states on $\clb_{t, \omega}$. i.e. a $\ast$-homomorphism. Indeed, if an extreme point $\zeta$ of $\cls$ can be written 
                 as $q\zeta_1+(1-q)\zeta_2$, where $0<q<1$ and $\zeta_1, \zeta_2$ are states on $\clb_{t, \omega}$, we have ${\rm ev}_p=q\zeta_1^\prime+(1-q)\zeta_2^\prime$, 
                  where $\zeta_i^\prime$ denotes the restriction of $\zeta_i$ to $C(M) \ot 1$. As ${\rm ev}_p$ is a pure state of $C(M) \ot 1$, this implies $\zeta^\prime_i={\rm ev}_p$ for 
                  $i=1,2$, i.e. $\zeta_i \in \cls$. Then, by  the extremality of $\zeta$ in $\cls$, $\zeta_i=\zeta$ for $i=1,2$. Hence $\zeta$ is a pure state of $\clb_{t, \omega}$, i.e. 
                   $\ast$-homomorphism and we have $\zeta(\Pi_t(f))=\phi(\zeta(Y_1(t)), \ldots, \zeta(Y_n(t)))$. 
                   
                   Now, recall from Corollary \ref{cor_state} that $({\rm id} \ot h)(\clb_{t, \omega})\subseteq C(M)$, so $\eta:=({\rm ev}_p \ot h)$ is a well-defined state
                    on $\clb_{t, \omega}$ and it is also 
                    an element of $\cls$. 
                    Moreover, as $f(p)=0$, (\ref{loc_formula}) 
                    implies the following :
                    $$ \hat{\cll}(f)(p)=
                    =\eta(U^{-1}\tilde{\cll}(\alpha(f))U)=
                   \eta \left( \lim_{\epsilon \raro 0+}  \frac{\IE ( \Pi_{\tau_\epsilon}(f))}{\IE(\tau_\epsilon)}\right)
                    =\lim_{\epsilon \raro 0+}  \frac{\IE ( \eta(\Pi_{\tau_\epsilon}(f)))}{\IE(\tau_\epsilon)}
                    .$$ 
              We claim that $$ \zeta( \Pi_{\tau_\epsilon}(f))=0~\forall \zeta \in \cls,$$ for all sufficiently small $\epsilon$.   It is enough to prove  it when $\zeta$ is an extreme point, i.e. $\ast$-homomorphism. For any such extremal state $\zeta$, 
                we have  
                $|\zeta(Y_i(\tau_\epsilon)-x_i \ot 1)|\leq \epsilon$ $\forall i$ by the continuity of the Brownian flow. As $\zeta(x_i \ot 1)=p_i$ by definition of $\cls$,  
                  the tuple 
               $(\zeta(Y_1({\tau_\epsilon})), \ldots, \zeta(Y_n({\tau_\epsilon})))  \in \IR^n $  is contained in an $n$-cube of side-length $\epsilon$ around 
                $(p_1, \ldots, p_n)$. Moreover, as $\zeta \circ \Pi_{\tau_\epsilon}$ is a character 
                of $C(M)$,  there is some point $v=(v_1, \ldots, v_n) \in M \subset \IR^n$ such that 
                 $ \zeta \circ \Pi_{\tau_\epsilon}(f)=f(v)$ for all $f \in C(M)$. In particular,  $(\zeta(Y_1({\tau_\epsilon})), \ldots, \zeta(Y_n({\tau_\epsilon})))=(v_1, \ldots, v_n) \in M$.
                                   Thus, $\zeta(\Pi_{\tau_\epsilon}(f))=\phi(\zeta(Y_1({\tau_\epsilon})), \ldots, \zeta(Y_n({\tau_\epsilon})))=0$ for all 
                  $\epsilon < \epsilon_0$, proving our claim. In particular, we have  $\eta(\Pi_{\tau_\epsilon}(f))=0$ for 
                  all sufficiently small $\epsilon$, hence $\hat{\cll}(f)(p)=0$.\vspace{2mm}\\

              {\it Step 3: Non-degeneracy on $\Omega^1(\clc_0)$.}\vspace{2mm}\\
     Let $f_1, \ldots, f_k,~g_1, \ldots, g_k$ ($k \geq 1$) be smooth functions in $\clc_0$ such that $\sum_{i,j=1}^k \overline{f_i}f_jk_{\hat{\cll}}(g_i,g_j)=0$. It follows from the proof of {\it Step 1} that 
     $$\sum_{i,j=1}^k \overline{f_i}f_jk_{\hat{\cll}}(g_i,g_j)
     =({\rm id} \ot h)\left( U^{-1} << \Omega, \Omega>>U\right),$$ where $\Omega=\sum_{i=1}^k d\alpha(g_i)\alpha(f_i).$ As $h$ is faithful, $({\rm id} \ot h)$ is so and therefore,  we get 
     $U^{-1} << \Omega, \Omega>>U=0$, hence $<<\Omega, \Omega>>=0$ which implies $\Omega=0$ as $<<\cdot, \cdot>>$ comes from a Riemannian structure and hence is  an inner product. Fix any $x \in M$ and a set of local coordinates $x_1, \ldots x_m$ around $x$. Then, $\Omega(x)=0$ implies 
     $\sum_{i=1}^k (\frac{\partial}{\partial x_j}g_{i(0)})(x)f_{i(0)}(x)\ot g_{i(1)}f_{i(1)}=0$ $\forall j$. Applying the counit $\epsilon$,  we get $$\sum_i \frac{\partial}{\partial x_j}
     (g_{i(0)}\epsilon(g_{i(1)}))(x)f_{i(0)}(x)
     \epsilon(f_{i(1)})=\sum_i (\frac{\partial}{\partial x_j}g_i)(x)f_i(x)=\sum_i f_i(x) (\frac{\partial}{\partial x_j}g_i)(x)=0$$
     $\forall j=1,\ldots, m,$ which means $\sum_i f_idg_i|_x=0$. As $x$ is arbitrary, it follows that $\sum_i f_i dg_i=0$.

                  \qed\\

                \bcrlre
                \label{metric_pres}
                Any smooth action on a compact Riemannian manifold preserves some Riemannian metric on $M$.
                \ecrlre
              {\it Proof:}\\ By Proposition \ref{nondeg_loc} and Theorem \ref{pre_main},   $\hat{\cll}$  induces a  $C(M)$-valued sesquilinear form given by 
              $<<df,dg>>^\prime=\hat{\cll}(\overline{f}g)-\hat{\cll}(\overline{f})g
                -\overline{f}\hat{\cll}(g) $ for all  smooth functions  $f,g$. We claim that $\alpha$ preserves this  sesquilinear form. It is also clear from the definition of $<<\cdot, \cdot>>^\prime$ that it 
                 is Fr\'echet continuous in the sense discussed in Subsection 3.1. 
                As $\alpha$ is smooth, it is enough to prove $<<df_{(0)},dg_{(0)}>>^\prime\ot f_{(1)}g_{(1)}=
                \alpha(<<df,dg>>^\prime)$ for all  $f,g \in \clc_0$. 
                 For this, it is enough to prove that $(\hat{\cll} \ot {\rm id})(\alpha(f))=\alpha(\hat{\cll}(f))$ for all $f \in \clc_0$ (hence 
                  for all $f \in C^\infty(M)$ by appropriate continuity of $\alpha$ and $\hat{\cll}$). Once we prove this, the argument of Lemma 4.3 
                   of \cite{gafa} can be applied verbatim. There $\cll$ is the Laplacian of a Riemannian structure but that has no role in the proof; 
                    the algebraic calculation requires only that $\cll$ commutes with $\alpha$. 
                    
                    Now, $\hat{\cll}(f)=\Psi(G)$ where $G=\cll(f_{(0)})\ot f_{(1)}$ 
                    and we have  the following by (\ref{psi_inv}) 
                     \bean \lefteqn {(\hat{\cll} \ot {\rm id} )(\alpha(f))}\\
                     &=& (\Psi \ot {\rm id})(\cll(f_{(0)(0)}) \ot f_{(0)(1)} \ot f_{(1)})\\
                     &=& (\Psi \ot {\rm id})(\cll(f_{(0)})\ot f_{(1)(1)} \ot f_{(1)(2)})\\
                     &=& (\Psi \ot {\rm id})({\rm id} \ot \Delta)(\cll(f_{(0)})\ot f_{(1)})\\
                     &=& \alpha(\Psi(\cll(f_{(0)})\ot f_{(1)}))\\
                     &=& \alpha(\hat{\cll}(f)).\eean
    
    We can adapt the proof of $(ii) \Rightarrow (i)$ of Theorem 3.5 of \cite{gafa} to conclude that (i) of that theorem holds, i.e. $d\alpha(g)(x) \alpha(f)(x)=\alpha(f)(x) d\alpha(g)(x)$ for any $f,g \in \clc_0,x \in M$ in the notation of \cite{gafa}. 
    Indeed, as in the proof of Theorem 3.5 of \cite{gafa}, we can take $F=\alpha(f) d\alpha(g)-d\alpha(g)\alpha(f)\in \Omega^1(\clc_0)$ where $f,g \in \clc_0$, and observe that the (purely algebraic) proof of the fact $<<F,F>>^\prime=0$ using the $<<\cdot, \cdot>>^\prime$-preservation of $\alpha$ goes through verbatim even if the sesquilinear form  does not come from a Riemannian structure. Moreover, as $<<\cdot,\cdot>>^\prime$ is an inner product on $\Omega^1(\clc_0)$,
     we get $F=0$. By Fr\'echet continuity of $\alpha$, we extend $d\alpha(g)(x) \alpha(f)(x)=\alpha(f)(x) d\alpha(g)(x)$ to all $f,g \in C^\infty(M)$.

    But then, Theorem 3.5 of \cite{gafa} gives us a Riemannian structure such that the corresponding inner product will be preserved by $\alpha$ as well. 
       \qed\\
       
   We have already observed the following in  the  proof of the above Corollary \ref{metric_pres}, which can be called `commutativity of partial derivatives up to the first order':
   \bcrlre
   \label{first_comm}
   For any point $x\in M$ and local coordinates $(x_1, \ldots, x_m)$ around $x$, the algebra $\clq_x$ generated by $\alpha(f)(x), \frac{\partial}{\partial x_i}\alpha(g)(x)$, where 
   $f,g \in C^\infty(M)$ and $i=1, \ldots,m$, is commutative.
   \ecrlre

   
   \subsection{Proof of the conjecture}
   Let $\alpha$ be a smooth action as in the previous subsection.  We have already seen commutativity of partial derivatives up to the first order. We want to prove similar commutativity for higher order partial derivatives. This involves  lift to the cotangent bundle.
   \blmma
   \label{higher_comm}
   For any point $x\in M$ and local coordinates $(x_1, \ldots, x_m)$ around $x$, the algebra  generated by $\alpha(f)(x), \frac{\partial}{\partial x_{i_1}}\ldots \frac{\partial}{\partial x_{i_k}}\alpha(g)(x)$, where $f,g \in C^\infty(M),$ $k \geq 1$ and $i_j \in \{1,\ldots, m\}$, is commutative.
   \elmma
   {\it Proof:}\\
      We need an analogue of Theorem 3.4 of \cite{gafa},  to lift the given action to a smooth action on the sphere bundle of the cotangent space. As the constructions and arguments in \cite{gafa} go through almost verbatim, we just sketch the main line of arguments very briefly. 
   
   First, we choose a Riemannian metric $<\cdot, \cdot>$  by Corollary \ref{metric_pres} 
 which is preserved by the action. Consider the compact smooth manifolds $S$  and $\tilde{S}$ given by:  $$S=\{ (x,\omega):~x\in M,~\omega \in T^*_xM~<\omega, \omega>_x=1\},$$ 
 $$\tilde{S}:=\{ (x,\omega):~x \in M,~ \omega \in T^*_x M,~\omega \ne 0\}.$$
 Let $\pi: S \raro M$ and $ \tilde{\pi}: \tilde{S} \raro M$ be the natural projection maps. In analogy with the construction of Subsection 3.3 of \cite{gafa}, we define $\theta_\xi \in C^\infty(S)$, where $\xi \in \Omega^1(C^\infty(M)),$ by $\theta_\xi(x,\omega)=<\omega,\xi(x)>_x$. For any local coordinate chart $(U,(x_1, \ldots, x_m))$ for $M$ and $U$-orthonormal one-forms $\omega^\prime_1, \ldots, \omega^\prime_m$ in the sense of \cite{gafa}, i.e. $\{ \omega^\prime_1(y),\ldots, \omega^\prime_m(y)\}$ is an orthonormal basis of $T^*_yM$ for all $y \in U$, we define $t_j^U=\theta_{\omega^\prime_j}.$ 
 Similarly, define $T_j^U \in C^\infty(S,\clq)$ by $$T_j^U(x,\omega)=
 <<\omega \ot 1,d\alpha_{(1)}(\omega^\prime_j)(x)>>,$$ where $d\alpha_{(1)}$ is the lift of $\alpha$ to the module of one-forms as in \cite{gafa} and $<<\cdot, \cdot>>$ denotes the $\clq$-valued inner product of $T^*_xM \ot \clq$. 
   Then, following the arguments of Subsection 3.3 of \cite{gafa}, we can prove that there exists a faithful smooth action $\beta$ (say) of $\clq$ on $S$. The action is determined by $$ \beta((f \circ \pi)t^U_j)=(\alpha(f)\circ \pi) T^U_j,$$ for any $f \in C_c^\infty(M)$ supported in a coordinate chart $U$. Equivalently, we have $\beta(\theta_{df})(x,\omega)=<<\omega \ot 1,d\alpha(f)(x)>>.$  Applying Corollary \ref{first_comm}
   to $\beta$, we conclude that for any $e \in S$, the algebra ($\clq^1_e$, say) generated by  $\{ (X \ot {\rm id})(\beta(F))(e),\beta(G)(e):~F,G \in C^\infty(S),~X \in \chi(S)\}$
   (where $\chi(S)$ denotes the set of smooth vector fields on $S$) is commutative. Let us extend $\beta$ further to 
    $\tilde{S}$.  Clearly, $\tilde{S}$ is diffeomorphic to $S \times \IR^{\times}$, where $\IR^\times=\IR \backslash \{ 0 \}$ and the diffeomorphism $\psi:\IR^{\times} \times
     S\raro \tilde{S} $ (say) is given by $\psi((x,\omega),r)=(x,r\omega)$. This induces the isomorphism $C_c(\tilde{S})
      \cong C_c(\IR^\times) \ot C(S)$. In what follows, we will interchangeably use the two equivalent descriptions of $\tilde{S}$ explained above, without explicitly mentioning the diffeomorphism $\psi$. 
      
      Define $\tilde{\beta}: C_c(\tilde{S}) \raro C_c(\tilde{S},\clq)$ by $\tilde{\beta}(\tilde{F})((x,\omega),r)=\beta(\tilde{F}_r)(x,\omega)$, where 
      $\tilde{F}_r \in C(S)$ is given by $\tilde{F}_r(x,\omega)=\tilde{F}((x,\omega),r)$.    
      
      From the definition it is clear that $\tilde{\beta}$ maps $C^\infty_c(\tilde{S})$ into $C_c^\infty(\tilde{S},\clq)$ and  
       for $\tilde{e}=(e,r) \in \tilde{S}$ ($e \in S$) any smooth vector field $Y$ on $\tilde{S}$ and smooth compactly supported  function $\tilde{F}$ on $\tilde{S} $,
        $(Y \ot {\rm id})(\tilde{\beta}(F))(\tilde{e})$ belongs to $\clq^1_e$. Indeed, it is enough to check this for $\tilde{F}$ of the form $\tilde{F}((x,\omega),r)=F(x,\omega) g(r)$, $F \in C^\infty(S)$, $g \in C_c^\infty(\IR^\times)$. For such $\tilde{F}$, we have $\tilde{\beta}(\tilde{F})((x,\omega),r)=\beta(F)(x,\omega)g(r).$ Moreover,  any smooth vector field $Y$ on $\tilde{S}$ can be written (locally) as $\phi_1 X +\phi_2 \frac{\partial}{\partial r}$ where $\phi_1,\phi_2$ are smooth functions on $\tilde{S}$ and $X$ is a vector field in the direction of $S$. Thus, $(Y \ot {\rm id})(\tilde{\beta}(\tilde{F}))(\tilde{e})=g(r)\phi_1(\tilde{e})(X \ot {\rm id})(\beta(F))(e)
        +g^\prime(r)\phi_2(\tilde{e})\beta(F)(e) \in \clq^1_e$.
        
         For a set of local coordinates $(U,(x_1, \ldots, x_m))$ for $M$ and $U$-orthonormal one-forms $\omega^\prime_1, \ldots, \omega^\prime_m$ as before, define $\tilde{t}^U_j: \tilde{S} \raro \IR$ by $\tilde{t}_j^U(e,r)=rt^U_j(e).$ It is clear from the definition of $t^U_j$ that $\sum_{j=1}^m (\tilde{t}_j^U)^2=r^2$ for all $(e,r)$ with $\pi(e) \in U$. Moreover, $(x_1, \ldots, x_m, \tilde{t}^U_1, \ldots, \tilde{t}^U_m)$ is a set of local coordinates for $\tilde{S}$ on the neighbourhood ${\tilde{\pi}}^{-1}(U)\cong \pi^{-1} \times \IR^\times$. Let us write $y_i$ for $\tilde{t}^U_i$ and define 
         $\tilde{\theta}_\xi: \tilde{S} \raro \IR$ ($\xi \in \Omega^1(C^\infty(M))$)   by $\tilde{\theta}_\xi((x,\omega), r)=<r\omega, \xi>(x)=r\theta_\xi(x,\omega).$
         
         Fix a
        point $\tilde{e}_0=(x_1^0, \ldots, x_m^0, y_1^0, \ldots y_m^0)$ of $\tilde{S}$, a coordinate neighbourhood 
         $(U,(x_1, \ldots, x_m))$ of $M$ around $(x_1^0, \ldots, x_m^0)$, a set of $U$-orthonormal one-forms $\omega^\prime_1, \ldots, \omega^\prime_m$ and  functions $f \in C^\infty(M)$ and $g \in C_c^\infty(\IR^\times)$ such that $g=1$ in an open neighbourhood  of $r_0:=((y_1^0)^2+\ldots 
         +(y^0_m)^2)^{\frac{1}{2}}.$ Consider $\tilde{F}=g\tilde{\theta}_{df} \in C_c^\infty(\tilde{S})$. Clearly, on a sufficiently small 
          neighbourhood of $\tilde{e}_0$, we have (using the orthonormality of $\omega^\prime_j(x)$'s) 
          \bean \lefteqn{\tilde{\beta}(\tilde{F})((x,\omega),r)}\\
          &=& r<<\omega \ot 1, d\alpha(f)(x)>>=\sum_j <\omega, r\omega^\prime_j(x)>_x<<\omega^\prime_j \ot 1,
          d\alpha(f)>>(x)\\
          &=&\sum_j \tilde{t^U_j}(x,\omega,r)<<\omega^\prime_j \ot 1,
          d\alpha(f)>>(x).\eean
   In other words, writing $x=(x_1, \ldots, x_m), $ $y=(y_1, \ldots, y_m)$, we have 
   $$\tilde{\beta}(\tilde{F})(x,y)=\sum_j y_j \eta_j(x),$$ where $\eta_j(x)=<<\omega^\prime_j,
          d\alpha(f)>>(x).$ This implies, 
        $$\frac{\partial}{\partial x_i} \tilde{\beta}(\tilde{F})(x,y)=
        \sum_j y_j (\frac{\partial }{\partial x_i}\eta_j)(x).$$
        We have already seen that the left hand side of the above belongs to $\clq^1_e$. Therefore, fixing $x=x^0:=(x^0_1, \ldots, x^0_m)$ we get $\sum_j y_j C_j \in \clq^1_{(x^0,y)}$, where $C_j=(\frac{\partial }{\partial x_i}\eta_j)(x^0),$ for all $y$ in an open neighbourhood of $(y^0_1, \ldots, y^0_m)$ in $\IR^m$. 
        As $\clq^1_{(x^0,y)}$ is a commutative algebra by Corollary \ref{first_comm}, we have 
        $$\sum_{j\leq k=1}^m y_jy_k[C_j,C_k]=0.$$ Using the fact that  $\{ y_jy_k, j\leq k\}$ are linearly independent as $y_1, \ldots, y_m$  are the  coordinates for an $m$-dimensional open neighbourhood, we conclude $[C_j,C_k]=0$. Moreover,
        for any $\phi \in C^\infty(M)$, we have $\alpha(\phi)(x)=\beta(\phi \circ \pi)(x,\omega)$ for any $\omega$, hence $\alpha(\phi)(x) \in \clq^1_{(x,\omega)}$. It follows that    $\alpha(\phi)(x^0)$ commutes with $\sum_j y_j C_j$, and using the linear independence of the $y_i$, we get $[\alpha(\phi)(x^0),C_j]=0$ for all $j$. Similarly, $<<\omega^\prime_j \ot 1,d\alpha(\phi)>>(x^0)\in \clq^1_{(x^0,y)}$ and this helps us conclude the commutativity between $<<\omega^\prime_j \ot 1,d\alpha(\phi)>>(x^0)$ and $C_k$ for any $j,k=1, \ldots,m$. 
        In other words, we have proved the commutativity of  the algebra (say, $\clb^\alpha_2(x^0)$) generated by $\alpha(f_1)(x^0), <<\omega^\prime_j \ot 1,d\alpha(f_2)>>(x^0), \frac{\partial}{\partial x_i}<<\omega^\prime_k \ot 1,d\alpha(f_3)>>(x^0)$, $f_p \in C^\infty(M),i,j,k=1, \ldots,m$, $p=1,2$. For $\omega \in \Omega^1(C^\infty(M)),$ denote by  $X_\omega$ the vector field given by $X_\omega(f)=<<\omega, df>>$ as in \cite{gafa}.  
        Writing $\frac{\partial}{\partial x_i}$ in terms of $X_{\omega^\prime_j}$'s (see also the arguments in the beginning of Theorem 4.6 of \cite{gafa}) 
         we can see that $\clb^\alpha_2(x^0)$ is the  same as the algebra generated by 
         $\alpha(f)(x^0), \frac{\partial}{\partial x_i}\alpha(g)(x^0) \frac{\partial}{\partial x_{j}} \frac{\partial}{\partial x_{k}}\alpha(\phi)(x^0)$, where $f,g, \phi \in C^\infty(M),$ $i,j,k=1,\ldots,m$. 
         
         We can go on like this and set up an induction hypothesis that the algebra $\clb^\alpha_l(x)$ (say) generated by $\alpha(f)(x), \frac{\partial}{\partial x_{i_1}}\ldots \frac{\partial}{\partial x_{i_k}}\alpha(g)(x)$, where $f,g \in C^\infty(M),$ $1\leq k \leq l,$ is commutative for any smooth action $\alpha$ on a compact smooth manifold $M$. Using the induction hypothesis (for $l$) for $\beta$ on  $S$, we see that $\frac{\partial}{\partial x_{i_1}}\ldots \frac{\partial}{\partial x_{i_l}}\tilde{\beta}(\tilde{F})(x,y)$ belongs to 
          a commutative algebra $\clb^\beta_l(e)$. Proceeding as before, 
          we conclude the commutativity of $\clb^\alpha_{l+1}$. 
        \qed\\
        
   \bthm
       \label{main}
       Let $\alpha$ be a smooth faithful action of a CQG $\clq$ on a compact connected smooth manifold $M$. Then $\clq$ must be 
        classical, i.e. isomorphic with $C(G)$ for  a compact group $G$ acting smoothly on $M$. 
       \ethm
{\it Proof}\\ Note that in the proof of Theorem 5.3 of \cite{gafa}, the isometry condition, i.e. commutation with the Laplacian, was used only to get commutativity of all order partial derivatives of the action. However, we have already proved this commutativity in Lemma \ref{higher_comm}. This allows the proof of Theorem 5.3 of \cite{gafa} to be carried through more or less verbatim. Let us sketch it briefly. 

Given the smooth action $\alpha$ of $\clq$ on $M$, we choose a Riemannian metric by Corollary \ref{metric_pres} 
 which is preserved by the action. 
 This implies the commutativity of $\clq_x$. Using this,  we can proceed along the lines of \cite{gafa} to lift the given action to $O(M)$.  
 Now, by Lemma \ref{higher_comm}, we do have the commutativity of partial derivatives of all orders for   the lifted action $\Phi$ needed in steps (i) and (iv) of  the proof of Theorem 5.3 of \cite{gafa}  and 
 the rest of the arguments of Theorem 5.3 of \cite{gafa} will go through.\qed\\
\brmrk 
 Observe that in the proof of Lemma 5.1 of \cite{gafa}, only commutativity of partial derivatives up to the second order is necessary. This means it is actually sufficient to state and prove Lemma \ref{higher_comm} for commutativity up to the second order. 
 \ermrk
As an application, we can generalize the results obtained by Chrivasitu in \cite{chirvasitu} for some other class of   Riemannian manifolds. More precisely, 
\bcrlre
Let $M$ be any compact connected Riemannian manifold so that the metric space $(M,d)$ (where $d$ is the Riemannian geodesic distance) 
 satisfies the hypotheses of Corollary 4.9 of \cite{Metric}. Then the quantum isometry group $QISO(M,d)$ in the sense of \cite{Metric}  coincides with $C(ISO(M,d)).$
\ecrlre
{\it Proof:}\\ It follows from the proof  of existence of  $QISO(M,d)$ in \cite{Metric} that the action of $QISO(M,d)$ on $C(M)$ is affine w.r.t.  the coordinate functions coming from 
 any  embedding $M \subseteq \IR^N$ satisfying the conditions of Corollary 4.9 of \cite{Metric}. But this means that the action is smooth in our sense, hence by Theorem \ref{main} we complete the proof.\qed\\

 {\bf Acknowledgment}: The author would like to thank A. Chirvasitu for pointing out some corrections and other comments.   He is also grateful to P. Hajac for pointing out the reference \cite{free}. 
    
\end{document}